\documentclass[global]{svjour}
\usepackage{amssymb}
\usepackage{amsfonts}
\usepackage{amsmath}

\input{tcilatex}

\begin{document}

\title{Existence of equilibrium for multiobjective games in abstract convex
spaces}
\author{Monica Patriche}
\institute{University of Bucharest 
\email{monica.patriche@yahoo.com}%
}
\mail{University of Bucharest, Faculty of Mathematics and Computer Science, 14
Academiei Street, 010014 Bucharest, Romania}
\maketitle

\begin{abstract}
In this paper we use the minimax inequalities obtained by S. Park (2011) to
prove the existence of weighted Nash equilibria and Pareto Nash equilibria
of a multiobjective game defined on abstract convex spaces.
\end{abstract}


\noindent \textit{2000 Mathematics Subject Classification:}\texttt{\ }47H10,
55M20, 91B50.

\section{\textbf{INTRODUCTION}}

Recently, in [8], S. Park introduced a new concept of abstract convex space
and several classes of correspondences having the KKM property. With this
new concept, the KKM type correspondences were used to obtain coincidence
theorems, fixed point theorems and minimax inequalities. S. Park generalizes
and unifies most of important results in the KKM theory on G-convex spaces,
H-spaces, and convex spaces (for example, see [8]-[13]).

For the history of KKM literature, we must remind Ky Fan [3], who extended
the original KKM theorem to arbitrarily topological vector space. The
property of close-valuedness of related KKM correspondences was replaced
with more general concepts. In [7], Luc and al. have introduced the concept
of intersectionally closed-valued correspondences and in [13], S. Park has
obtained new KKM type theorems for this kind of KKM\ correspondences.

In this paper we use the minimax inequalities obtained by S. Park in [13] to
prove the existence of weighted Nash equilibria and Pareto Nash Equilibria
of a multiobjective game defined on abstract convex spaces. For the history
of minimax theorems, I also must remind the name of Ky Fan (see [4]). Among
the authors who studied the existence of Pareto equilibria in game theory
with vector payoffs, I emphasize S. Chebbi [2], W. K. Kim [5], W. K. Kim, X.
P. Ding [6], H. Yu [16], J. Yu, G. X.-Z Yuan [17], X. Z. Yuan, E. Tarafdar
[18]. A reference work is the paper of M. Zeleny [19]. The approaches of
above-mentioned authors deal with the Ky Fan minimax inequality,
quasi-equilibrium theorems or quasi-variational inequalities. We must
mention the papers of P. Borm, F. Megen, S. Tijs [1], who introduced the
concept of perfectness for multicriteria games and M. Voorneveld, S. Grahn,
M. Dufwenberg [14], who studied the existence of ideal equilibria. Ather
authors, as H. Yu (see [16]), obtained the existence of a solution of
multiobjective games by using new concepts of continuity and convexity.

The paper is organised as follows: In section 2, some notation,
terminological convention, basic definitions and results about abstract
convex spaces and minimax inequalities are given. Section 3 introduces the
model, that is, a multiobjective game defined on an abstract convex space
and the concept of weight Nash equilibrium. {\normalsize Section 4 contains
existence results for weight Nash equilibrium and Pareto Nash equilibrium.}

\section{\textbf{ABSTRACT CONVEX SPACES AND MINIMAX INEQUALITIES}}

Let $A$ be a subset of a topological space $X$. $2^{A}$ denotes the family
of all subsets of $A$. $\overline{A}$ denotes the closure of $A$ in $X$ and
int$A$ denotes the interiorof $A$. If $A$ is a subset of a vector space, co$%
A $ denotes the convex hull of $A$. If $F$, $G:$ $X\rightarrow 2^{Y}$ are
correspondences, then co$G$, cl $G$, $G\cap F$ $:$ $X\rightarrow 2^{Y}$ are
correspondences defined by $($co$G)(x)=$co$G(x)$, $($cl$G)(x)=$cl$G(x)$ and $%
(G\cap F)(x)=G(x)\cap F(x)$ for each $x\in X$, respectively. The graph of $%
F:X\rightarrow 2^{Y}$ is the set Gr$(F)=\{(x,y)\in X\times Y\mid y\in F(x)\}$
and $F^{-}:Y\rightarrow 2^{X}$ is defined by $F^{-}(y)=\{x\in X:y\in F(x)\}$
for $y\in Y.$ Let $\tciFourier (A)$ be the set of all nonempty finite
subsets of a set $A.\medskip $

For the reader's convenience, we review a few basic definitions and results
from abstract convex spaces.

\textit{Definition 1 [13].} Let $X$ be a topological space, $D$ be a
nonempty set and let $\Gamma :\tciFourier (D)\rightarrow 2^{X}$ be a
correspondence with nonempty values $\Gamma _{A}=\Gamma (A)$ for $A\in
\tciFourier (D).$ The family $(X,D;\Gamma )$ is called an \textit{abstract
convex space.}

\textit{Definition 2 [13]. }For a nonempty subset $D^{\prime }$ of $D$, we
define the $\Gamma $\textit{-convex hull of} $D^{\prime }$, denoted by co$%
_{\Gamma }D^{\prime }$, as

co$_{\Gamma }D^{\prime }$=$\cup \{\Gamma _{A}:A\in \tciFourier (D^{\prime
})\}\subset X.$

\textit{Definition 3 [13]. }Given an abstract convex space ($X,D,\Gamma $),
a nonempty subset $Y$ of $X$ is called to be a $\Gamma $\textit{-convex} 
\textit{subset} of ($X,D,\Gamma $) relative to $D^{\prime }$ if for any $%
A\in \tciFourier $($D^{\prime }$), we have $\Gamma _{A}\subset Y$, that is,
co$_{\Gamma }D^{\prime }\subset Y.$

\textit{Definition 4 [13]. }When $D\subset X$ in ($X,D,\Gamma $), a subset $%
Y $ of $X$ is said to be $\Gamma $\textit{-convex }if co$_{\Gamma }(Y\cap
D)\subset Y;$ in other words, $Y$ is $\Gamma $\textit{-convex }relative to%
\textit{\ }$D^{\prime }=Y\cap D.$ In case $X=D,$ let $(X,\Gamma
)=(X,X,\Gamma ).$

\textit{Definition 5 [13]. }The abstract convex space ($X,D,\Gamma $) is
called \textit{compact} if $X$ is compact.

We have abstract convex subspaces as the following simple observation.

\begin{proposition}
For an abstract convex space ($X,D,\Gamma $) and a nonempty subset $%
D^{\prime }$ of $D$, let $Y$ be a $\Gamma $-convex subset of $X$ relative to 
$D^{\prime }$ and $\Gamma ^{\prime }:\tciFourier $($D^{\prime }$)$%
\rightarrow 2^{Y}$ a correspondence defined by $\ $
\end{proposition}

$\Gamma _{A}^{\prime }=\Gamma _{A}\subset X$ for $A\in \tciFourier $($%
D^{\prime }$).

Then ($Y,D^{\prime },\Gamma ^{\prime }$) itself is an abstract convex space
called a subspace relative to $D^{\prime }.\medskip $

The following result is known.

\begin{lemma}[12]
Let ($X_{i},D_{i},\Gamma _{i}$)$_{i\in I}$ be any family of abstract convex
spaces. Let $X=\tprod\nolimits_{i\in I}X_{i}$ be equipped with the product
topology and $D=\tprod\nolimits_{i\in I}D_{i}$. For each $i\in I$, let $\pi
_{i}:D\rightarrow D_{i}$ be the projection. For each $A\in \tciFourier $($D$%
), define $\Gamma (A)=\tprod\nolimits_{i\in I}\Gamma _{i}(\pi _{i}(A))$.
Then ($X,D,\Gamma $) is an abstract convex space. \medskip
\end{lemma}

\textit{Definition 6 [13]. }Let ($X,D,\Gamma $) be an abstract convex space.
Then $F:D\rightarrow 2^{X}$ is called a \textit{KKM correspondence} if it
satisfies $\Gamma _{A}\subset F(A):=\cup _{y\in A}F(y)$ for all $A\in
\tciFourier $($D$).\medskip\ 

\textit{Definition 7 [13]. }The \textit{partial KKM principle} for an
abstract convex space ($X,D,\Gamma $) is the statement that, for any
closed-valued KKM correspondence $F:D\rightarrow 2^{X}$, the family $%
\{F(z)\}_{z\in D}$ has the finite intersection property. The \textit{KKM\
principle} is the statement that the same property also holds for any
open-valued KKM correspondence.

An abstract convex space is called a \textit{KKM space} if it satisfies the
KKM principle.\medskip

\begin{proposition}
Let ($X,D,\Gamma $) be an abstract convex space and ($X,D^{\prime },\Gamma
^{\prime }$) a subspace. If ($X,D,\Gamma $) satisfies the partial KKM
principle, then so does ($X$, $D^{\prime },\Gamma ^{\prime }$).\medskip
\end{proposition}

Let ($X,D,\Gamma $) be an abstract convex space.

\textit{Definition 8 [13]. }The function $f:X\rightarrow \overline{\mathbb{R}%
}$ is said to be quasiconcave (resp. quasiconvex) if $\{x\in X:f(x)>r\}$
(resp., $\{x\in X:f(x)<r\}$ is $\Gamma $-convex for each $r\in \overline{%
\mathbb{R}}.\medskip $

In [7], Luc and al. have introduced the concept of intersectionally
closed-valued correspondences.

\textit{Definition 9. }Let $F:D\rightarrow 2^{X}$ be a correspondence.

(i) [7] $F$ is \textit{intersectionally closed-valued} if $\cap _{z\in D}%
\overline{F(z)}=\overline{\cap _{z\in D}F(z)};$

(ii) $F$ is \textit{transfer closed-valued} if $\cap _{z\in D}\overline{F(z)}%
=\cap _{z\in D}F(z);$

(iii) [7] $F$ is \textit{unionly open-valued} if Int$\cup _{z\in D}F(z)=\cup
_{z\in D}$Int$F(z);$

\ (iv) $F$ is \textit{transfer open-valued} if $\cup _{z\in D}F(z)=\cup
_{z\in D}$Int$F(z);$

Luc at al. [7] noted that (ii)$\Rightarrow $(i$).$

\begin{proposition}[7]
The correspondence $F$ is intersectionally closed-valued (resp. \textit{%
transfer closed-valued}) if only if its complement $F^{C}$ is \textit{%
unionly open-valued} (resp. \textit{transfer open-valued}).\medskip
\end{proposition}

\textit{Definition 10 [13].} Let$Y$ be a subset of $X.$

(i) $Y$ is said to be \textit{intersectionally closed} (resp. \textit{%
transfer closed}) if there is an intersectionally (resp., transfer)
closed-valued correspondence $F:D\rightarrow 2^{X}$ such that $Y=F(z)$ for
some $z\in D.$

(ii) $Y$ is said to be \textit{unionly open} (resp. \textit{transfer open})
if there is an unionly (resp., transfer) open-valued correspondence $%
F:D\rightarrow 2^{X}$ such that $Y=F(z)$ for some $z\in D.\medskip $

S. Park gives in [13] the concept of generally lower (resp. upper)
semicontinuous function.

\textit{Definition 11 [13]. }The function $f:D\times X\rightarrow \overline{%
\mathbb{R}}$ is said to be \textit{generally lower (resp. upper)
semicontinuous} (g.l.s.c.) (resp. g.u.s.c.) on $X$ whenever, for each $z\in
D,$ $\{y\in X:f(z,y)\leq r\}$ (resp., $\{y\in X:f(z,y)\geq r\})$ is
intersectionally closed for each $r\in \overline{\mathbb{R}}.\medskip $

The aim of this paper is to prove the existence of a weighted Nash
equilibrium for a multicriteria game defined in the framework of abstract
convex spaces. For our purpose, we need the following theorem (variant of
Theorem 6.3 in [13]).\medskip

\begin{theorem}
(Minimax inequality, [13]). Let ($X,D=X,\Gamma $) an abstract convex space
satisfying the partial KKM principle, $f,g:X\times X\rightarrow \overline{%
\mathbb{R}}$ extended real-valued functions and $\gamma \in \overline{%
\mathbb{R}}$ such that
\end{theorem}

\textit{(i) for each }$x\in X,$\textit{\ }$g(x,x)\leq \gamma ;$

\textit{(ii) for each }$y\in X,$\textit{\ }$F(y)=\{x\in X:f(x,y)\leq \gamma
\}$\textit{\ is intersectionally closed (respectiv, transfer closed);}

\textit{(iii) for each }$x\in X,$\textit{\ co}$_{\Gamma }\{y\in
X:f(x,y)>\gamma \}\subset \{y\in X:g(x,y)>\gamma \};$

\textit{(iv) the correspondence }$F:X\rightarrow 2^{X}$\textit{\ satisfies
the following condition:}

\textit{\ \ \ \ there exists a nonempty compact subset }$K$\textit{\ of }$X$%
\textit{\ such that either}

\textit{\ \ \ \ (a) }$K\supset \cap \{\overline{F(y)}:y\in M\}$\textit{\ for
some }$M\in \tciFourier $\textit{(}$X$\textit{); or}

\textit{\ \ \ \ (b) for each }$N\in \tciFourier $\textit{(}$X$\textit{)}$,$%
\textit{\ there exists a compact }$\Gamma $\textit{-convex subset }$L_{N}$%
\textit{\ of }$X$\textit{\ relative to some }$X^{\prime }\subset X$\textit{\
such that }$N\subset X^{\prime }$\textit{\ and }$K\supset L_{N}\cap \cap
_{y\in X^{\prime }}\overline{F(y)}\neq \phi .$

\textit{Then}

\textit{1) there exists a }$x_{0}\in X$\textit{\ }$($\textit{resp., }$%
x_{0}\in K)$\textit{\ such that }$f(x_{0},y)\leq \gamma $\textit{\ for all }$%
y\in X;$

\textit{2) if }$\gamma :=\sup_{x\in X}g(x,x)$\textit{, then we have}

\textit{\ \ \ inf}$_{x\in X}\sup_{y\in X}f(x,y)\leq \sup_{x\in
X}g(x,x).\medskip $

For the case when $X=D$ (we are concerned with compact abstract spaces ($%
X,\Gamma $) satisfying the partial KKM principle), we have the following
variants of the corollaries stated in [13].

\begin{corollary}[13]
Let $f,$ $g:X\times X\rightarrow \mathbb{R}$ be real-valued functions and $%
\gamma \in \mathbb{R}$ such that
\end{corollary}

\textit{(i) for each }$x,y\in X,$\textit{\ }$f(x,y)\leq g(x,y)$\textit{\ and 
}$g(x,x)\leq \gamma ;$

\textit{(ii) for each }$y\in X,$\textit{\ }$\{x\in X:f(x,y)>\gamma \}$%
\textit{\ is unionly open in }$X$\textit{;}

\textit{(iii) for each }$x\in X,$\textit{\ }$\{y\in X:g(x,y)>\gamma \}$%
\textit{\ is }$\Gamma $\textit{-convex on }$X;$

\textit{Then}

\textit{1) there exists a }$x_{0}\in X$\textit{\ such that }$f(x_{0},y)\leq
\gamma $\textit{\ for all }$y\in X;$

\textit{2) if }$\gamma :=\sup_{x\in X}g(x,x)$\textit{, then we have}

\textit{\ \ \ inf}$_{x\in X}\sup_{y\in X}f(x,y)\leq \sup_{x\in
X}g(x,x).\medskip $

\begin{corollary}[13]
Let $f,$ $g:X\times X\rightarrow \mathbb{R}$ be functions such that
\end{corollary}

\textit{(i) for each }$x,y\in X,$\textit{\ }$f(x,y)\leq g(x,y)$ and $%
g(x,x)\leq \gamma ;$

\textit{(ii) for each }$y\in X,$\textit{\ }$f(\cdot ,y)$\textit{\ is g.l.s.c
on }$X$\textit{;}

\textit{(iii) for each }$x\in X,$\textit{\ }$f(x,\cdot )$\textit{\ is
quasiconcave on }$X;$

\textit{Then we have}

\textit{\ \ \ inf}$_{x\in X}\sup_{y\in X}f(x,y)\leq \sup_{x\in
X}g(x,x).\medskip $

\section{\textbf{MULTIOBJECTIVE GAMES}}

Now we consider the multicriteria game (or multiobjective game) in its
strategic form. Let $I$ be a finite set of players and for each $i\in I,$
let $X_{i}$ be the set of strategies such that $X=\tprod\nolimits_{i\in
I}X_{i}$ and ($X_{i}$, $D_{i},\Gamma _{i}$ for each $i\in I$) is an abstract
convex space with $D_{i}\subset X_{i}$. Let $T^{i}:X\rightarrow 2^{\mathbb{R}%
^{k_{i}}}$, where $k_{i}\in \mathbb{N}$, which is called the payoff function
(or called multicriteria). From Lemma 1, we also have that $(X,D,\Gamma )$
is an abstract convex space, where $X=\tprod\nolimits_{i\in I}X_{i},$ $%
D=\tprod\nolimits_{i\in I}D_{i}$ and $\Gamma (A)=\tprod\nolimits_{i\in
I}\Gamma _{i}(\pi _{i}(A))$ for each $A\in \tciFourier $($D$).\medskip

\textit{Definition 12. }The family $G=((X_{i},D_{i},\Gamma
_{i}),T^{i})_{i\in I}$ is called multicriteria game.\medskip

If an action $x:=(x_{1},x_{2},...,x_{n})$ is played, each player $i$ is
trying to find his/her payoff function $%
T^{i}(x):=(T_{1}^{i}(x),...,T_{k_{i}}^{i}(x)),$ which consists of
noncommensurable outcomes. We assume that each player is trying to minimize
his/her own payoff according with his/her preferences.

In order to introduce the equilibrium concepts of a multicriteria game, we
need several necessary notation.

\textit{Notation.} We shall denote by

$\mathbb{R}_{+}^{m}:=\{u=(u_{1},u_{2},...u_{m})\in \mathbb{R}^{m}:u_{j}\geq
0 $ $\forall j=1,2,...,m\}$ and

int$\mathbb{R}_{+}^{m}:=\{u=(u_{1},u_{2},...u_{m})\in \mathbb{R}^{m}:u_{j}>0$
$\forall j=1,2,...,m\}$

the non-negative othant of $\mathbb{R}^{m}$ and respective the non-empty
interior of $\mathbb{R}_{+}^{m}$ with the topology induced in terms of
convergence of vector with respect to the Euclidian metric.

\textit{Notation.} For each $i\in I,$ denote $X_{-i}:=\tprod\nolimits_{j\in
I\setminus \{i\}}X_{j}.$ If $x=(x_{1},x_{2},...,x_{n})\in X,$ we denote $%
x_{-i}=(x_{1},...,x_{i-1},x_{i+1},...,x_{n})\in X_{-i}.$ If $x_{i}\in X_{i}$
and $x_{-i}\in X_{-i}$, we shall use the notation $%
(x_{-i},x_{i})=(x_{1},...,x_{i-1},x_{i},x_{i+1},...,x_{n})=x\in X.$

\textit{Notation.} For each $u,v\in \mathbb{R}^{m}$, $u\cdot v$ denote the
standard Euclidian inner product.

Let $\widehat{x}=(\widehat{x}_{1},\widehat{x}_{2},...,\widehat{x}_{n})\in X.$
Now we have the following definitions.\medskip

\textit{Definition 13. }A strategy $\widehat{x}_{i}\in X_{i}$ of player $i$
is said to be a \textit{Pareto efficient strategy} (resp., a \textit{weak
Pareto efficient strategy}) with respect to $\widehat{x}\in X$ of the
multiobjective game $G=((X_{i},D_{i},\Gamma _{i}),T^{i})_{i\in I}$ if there
is no strategy $x_{i}\in X_{i}$ such that

$T^{i}(\widehat{x})-T^{i}(\widehat{x}_{-i},x_{i})\in \mathbb{R}%
_{+}^{k_{i}}\backslash \{0\}$ (resp., $T^{i}(\widehat{x})-T^{i}(\widehat{x}%
_{-i},x_{i})\in $int$\mathbb{R}_{+}^{k_{i}}\backslash \{0\}).\medskip $

\begin{remark}
Each Pareto equilibrium is a weak Pareto equilibrium, but the converse is
not always true.
\end{remark}

\textit{Definition 14. }A strategy $\widehat{x}\in X$ is said to be a 
\textit{Pareto equilibrium} (resp., a \textit{weak Pareto equilibrium}) of
the multiobjective game $G=((X_{i},D_{i},\Gamma _{i}),T^{i})_{i\in I}$ if
for each player $i\in I$, $\widehat{x}_{i}\in X_{i}$ is a Pareto efficient
strategy (resp., a weak Pareto efficient strategy) with respect to $\widehat{%
x}.\medskip $

\textit{Definition 15. }A strategy $\widehat{x}\in X$ is said to be a 
\textit{weighted Nash} \textit{equilibrium} with respect to the weighted
vector $W=(W_{i})_{i\in I}$ with $W_{i}=(W_{i,1},W_{i,2},...,W_{i,k_{i}})\in 
\mathbb{R}_{+}^{k_{i}}$ of the multiobjective game $G=((X_{i},D_{i},\Gamma
_{i}),T^{i})_{i\in I}$ if for each player $i\in I$, we have

(i) $W_{i}\in \mathbb{R}_{+}^{k_{i}}\backslash \{0\};$

(ii) $W_{i}\cdot T^{i}(\widehat{x})\leq W_{i}\cdot T^{i}(\widehat{x}%
_{-i},x_{i}),$ $\forall x_{i}\in X_{i},$where $\cdot $ denotes the inner
product in $\mathbb{R}^{k_{i}}.$

\begin{remark}
In particular, if $W_{i}\in \mathbb{R}_{+}^{k_{i}}$ with $%
\tsum\nolimits_{j=1}^{k_{i}}W_{i,j}=1$ for each $i\in I,$ then the strategy $%
\widehat{x}\in X$ is said to be a normalized weighted Nash equilibrium with
respect to $W.\medskip $
\end{remark}

\section{\textbf{EXISTENCE OF WEIGHTED NASH EQUILIBRIUM AND\ PARETO\ NASH
EQUILIBRIUM}}

Now, as an application of Theorem 1, we have the following existence theorem
of weighted Nash equilibria for multiobjective games.

\begin{theorem}
Let $I$ be a finite set of indices, let ($X_{i},D_{i}=X_{i},\Gamma _{i}$)$%
_{i\in I}$ be any finite family of abstract convex spaces such that the
product space ($X,\Gamma $) satisfies the partial KKM principle. If there is
a weighted vector $W=(W_{1},W_{2},...,W_{n})$ with $W_{i}\in \mathbb{R}%
_{+}^{k_{i}}\backslash \{0\}$ such that the followings are satisfied:
\end{theorem}

\textit{(i) for each }$y\in X,$\textit{\ }$F(y)=\{x\in
X:\tsum\nolimits_{i=1}^{n}W_{i}\cdot
(T^{i}(x_{-i},x_{i})-T^{i}(x_{-i},y_{i}))\leq 0\}$\textit{\ is
intersectionally closed (respectiv, transfer closed);}

\textit{(ii) there exists }$g:X\times X\rightarrow \overline{\mathbb{R}}$%
\textit{\ extended real-valued function such that for each }$x\in X,$\textit{%
\ }$g(x,x)\leq 0$\textit{\ and for each }$x\in X,$\textit{\ co}$_{\Gamma
}\{y\in X:\tsum\nolimits_{i=1}^{n}W_{i}\cdot
(T^{i}(x_{-i},x_{i})-T^{i}(x_{-i},y_{i}))>0\}\subset \{y\in X:g(x,y)>0\};$

\textit{(iii) the correspondence }$F:X\rightarrow 2^{X}$\textit{\ satisfies
the following condition:}

\textit{\ \ \ \ there exists a nonempty compact subset }$K$\textit{\ of }$X$%
\textit{\ such that either}

\textit{\ \ \ \ (a) }$K\supset \cap \{\overline{F(y)}:y\in M\}$\textit{\ for
some }$M\in \tciFourier $\textit{(}$X$\textit{); or}

\textit{\ \ \ \ (b) for each }$N\in \tciFourier $\textit{(}$X$\textit{)}$,$%
\textit{\ there exists a compact }$\Gamma $\textit{-convex subset }$L_{N}$%
\textit{\ of }$X$\textit{\ relative to some }$X^{\prime }\subset X$\textit{\
such that }$N\subset X^{\prime }$\textit{\ and }$K\supset L_{N}\cap \cap
_{y\in x^{\prime }}\overline{F(y)}\neq \phi ;$

\textit{then there exists }$\widehat{x}\in K$\textit{\ such that }$\widehat{x%
}$\textit{\ is a weighted Nash equilibria of the game }$G=((X_{i},\Gamma
_{i}),T^{i})_{i\in I}$\textit{\ with respect to }$W.\medskip $

\textit{Proof.} Define the function $f:X\times X\rightarrow \mathbb{R}$ by $%
f(x,y)=\tsum\nolimits_{i=1}^{n}W_{i}\cdot
(T^{i}(x_{-i},x_{i})-T^{i}(x_{-i},y_{i})),$ $(x,y)\in X\times X.$ By Theorem
1, we have that inf$_{x\in X}\sup_{y\in X}f(x,y)\leq \sup_{x\in X}g(x,x)=0.$
It follows that there exists an $\widehat{x}\in K$ such that $f(\widehat{x}%
,y)\leq 0$ for any $y\in X.$ That is $\tsum\nolimits_{i=1}^{n}W_{i}\cdot
(T^{i}(\widehat{x}_{-i},\widehat{x}_{i})-T^{i}(\widehat{x}_{-i},y_{i})\leq 0$
for any $y\in X.$ For any given $i\in I$ and any given $y_{i}\in X_{i},$ let 
$y=(\widehat{x}_{-i},y_{i}).$ Then we have

$W_{i}\cdot (T^{i}(\widehat{x}_{-i},\widehat{x}_{i})-T^{i}(\widehat{x}%
_{-i},y_{i}))=$

$=\tsum\nolimits_{j=1}^{n}W_{j}\cdot (T^{j}(\widehat{x}_{-i},\widehat{x}%
_{i})-T^{i}(\widehat{x}_{-i},y_{i}))-\tsum\nolimits_{j\neq i}W_{j}\cdot
(T^{j}(\widehat{x}_{-i},\widehat{x}_{i})-T^{i}(\widehat{x}_{-i},y_{i}))$

$=\tsum\nolimits_{j=1}^{n}W_{j}\cdot (T^{j}(\widehat{x}_{-i},\widehat{x}%
_{i})-T^{i}(\widehat{x}_{-i},y_{i}))\leq 0.$

Therefore, we have $W_{i}\cdot (T^{i}(\widehat{x}_{-i},\widehat{x}%
_{i})-T^{i}(\widehat{x}_{-i},y_{i}))$ $\leq 0$ for each $i\in I$ and $%
y_{i}\in X_{i},$ that is $\widehat{x}\in K$ is a weighted Nash equilibrium
of the game $G$ with respect to $W.\medskip $

We obtain the following corollaries for the compact games when $X=D$.

\begin{corollary}
Let $I$ be a finite set of indices, let ($X_{i},\Gamma _{i}$)$_{i\in I}$ be
any finite family of abstract convex spaces such that the product space ($%
X,\Gamma $) satisfies the partial KKM principle. If there is a weighted
vector $W=(W_{1},W_{2},...,W_{n})$ with $W_{i}\in \mathbb{R}%
_{+}^{k_{i}}\backslash \{0\}$ such that the followings are satisfied:
\end{corollary}

\textit{(i) there exists }$g:X\times X\rightarrow R$\textit{\ such that for
each }$x,y\in X,$\textit{\ }$\tsum\nolimits_{i=1}^{n}W_{i}\cdot
(T^{i}(x_{-i},x_{i})-T^{i}(x_{-i},y_{i}))\leq g(x,y)$\textit{\ and }$%
g(x,x)\leq 0;$

\textit{(ii) for each }$y\in X,$\textit{\ }$\{x\in
X:\tsum\nolimits_{i=1}^{n}W_{i}\cdot
(T^{i}(x_{-i},x_{i})-T^{i}(x_{-i},y_{i}))>0\}$\textit{\ is unionly open in }$%
X$\textit{;}

\textit{(iii) for each }$x\in X,$\textit{\ }$\{y\in X:g(x,y)>0\}$\textit{\
is }$\Gamma $\textit{-convex on }$X;$

\textit{then there exists }$\widehat{x}\in X$\textit{\ such that }$\widehat{x%
}$\textit{\ is a weighted Nash equilibria of the game }$G=((X_{i},\Gamma
_{i}),T^{i})_{i\in I}$\textit{\ with respect to }$W.\medskip $

\begin{corollary}
Let $I$ be a finite set of indices, let ($X_{i},\Gamma _{i}$)$_{i\in I}$ be
any finite family of abstract convex spaces such that the product space ($%
X,\Gamma $) satisfies the partial KKM principle. If there is a weighted
vector $W=(W_{1},W_{2},...,W_{n})$ with $W_{i}\in \mathbb{R}%
_{+}^{k_{i}}\backslash \{0\}$ such that the followings are satisfied:
\end{corollary}

\textit{(i) there exists }$g:X\times X\rightarrow R$\textit{\ such that for
each }$x,y\in X,$\textit{\ }$\tsum\nolimits_{i=1}^{n}W_{i}\cdot
(T^{i}(x_{-i},x_{i})-T^{i}(x_{-i},y_{i}))\leq g(x,y);$

\textit{(ii) for each fixed }$y\in X,$\textit{\ the function }$x\rightarrow
\tsum\nolimits_{i=1}^{n}W_{i}\cdot (T^{i}(x_{-i},x_{i})-T^{i}(x_{-i},y_{i}))$%
\textit{\ is g.l.s.c on }$X$\textit{;}

\textit{(iii) for each fixed }$x\in X,$\textit{\ the function }$y\rightarrow
\tsum\nolimits_{i=1}^{n}W_{i}\cdot (T^{i}(x_{-i},x_{i})-T^{i}(x_{-i},y_{i}))$%
\textit{\ is quasiconcave on }$X$\textit{;}

\textit{then there exists }$\widehat{x}\in X$\textit{\ such that }$\widehat{x%
}$\textit{\ is a weighted Nash equilibria of the game }$G=((X_{i},\Gamma
_{i}),T^{i})_{i\in I}$\textit{\ with respect to }$W.\medskip $

In order to prove an existence theorem of Pareto equilibria for
multiobjective games, we need the following lemma.

\begin{lemma}[15]
Each normalized weighted Nash equilibrium $\widehat{x}\in X$ with a weight $%
W=(W_{1},W_{2},...,W_{n})$ with $W_{i}\in \mathbb{R}_{+}^{k_{i}}\backslash
\{0\}$ (resp., $W_{i}\in $int$\mathbb{R}_{+}^{k_{i}}\backslash \{0\})$ and $%
\tsum\nolimits_{j=1}^{k_{i}}W_{i,j}=1$ for each $i\in I,$ for a
multiobjective game $G=(X_{i},T^{i})_{i\in I}$ is a weak Pareto equilibrium
(resp. a Pareto equilibrium) of the game $G.$

\begin{remark}
The conclusion of Lemma2 still holds if $\widehat{x}\in X$ is a weighted
Nash equilibrium with a weight $W=(W_{1},W_{2},...,W_{n})$, $W_{i}\in 
\mathbb{R}_{+}^{k_{i}}\backslash \{0\}$ for $i\in I$(resp., $W_{i}\in $int$%
\mathbb{R}_{+}^{k_{i}}\backslash \{0\}$ for $i\in I)$ of the game $G.$
\end{remark}

\begin{remark}
A Pareto equilibrium of $G$ is not necessarily a weighted Nash equilibrium
of the game $G.\medskip $
\end{remark}
\end{lemma}

\begin{theorem}
Let $I$ be a finite set of indices, let ($X_{i},D_{i}=X_{i},\Gamma _{i}$)$%
_{i\in I}$ be any finite family of abstract convex spaces such that the
product space ($X,\Gamma $) satisfies the partial KKM principle. If there is
a weighted vector $W=(W_{1},W_{2},...,W_{n})$ with $W_{i}\in \mathbb{R}%
_{+}^{k_{i}}\backslash \{0\}$ such that the followings are satisfied:
\end{theorem}

\textit{(i) for each }$y\in X,$\textit{\ }$F(y)=\{x\in
X:\tsum\nolimits_{i=1}^{n}W_{i}\cdot
(T^{i}(x_{-i},x_{i})-T^{i}(x_{-i},y_{i}))\leq 0\}$\textit{\ is
intersectionally closed (respectiv, transfer closed);}

\textit{(ii) there exists }$g:X\times X\rightarrow \overline{\mathbb{R}}$%
\textit{\ extended real-valued function such that for each }$x\in X,$\textit{%
\ }$g(x,x)\leq 0$\textit{\ and for each }$x\in X,$\textit{\ co}$_{\Gamma
}\{y\in X:\tsum\nolimits_{i=1}^{n}W_{i}\cdot
(T^{i}(x_{-i},x_{i})-T^{i}(x_{-i},y_{i}))>0\}\subset \{y\in X:g(x,y)>0\};$

\textit{(iii) the correspondence }$F:X\rightarrow 2^{X}$\textit{\ satisfies
the following condition:}

\textit{\ \ \ \ there exists a nonempty compact subset }$K$\textit{\ of }$X$%
\textit{\ such that either}

\textit{\ \ \ \ (a) }$K\supset \cap \{\overline{F(y)}:y\in M\}$\textit{\ for
some }$M\in \tciFourier $\textit{(}$X$\textit{); or}

\textit{\ \ \ \ (b) for each }$N\in \tciFourier $\textit{(}$X$\textit{)}$,$%
\textit{\ there exists a compact }$\Gamma $\textit{-convex subset }$L_{N}$%
\textit{\ of }$X$\textit{\ relative to some }$X^{\prime }\subset X$\textit{\
such that }$N\subset X^{\prime }$\textit{\ and }$K\supset L_{N}\cap \cap
_{y\in x^{\prime }}\overline{F(y)}\neq \phi ;$

\textit{then there exists }$\widehat{x}\in K$\textit{\ such that }$\widehat{x%
}$\textit{\ is a weak Pareto equilibrium of the game }$%
G=((X_{i},D_{i}=X_{i},\Gamma _{i}),T^{i})_{i\in I}.$\textit{\ In addition,
if }$W=(W_{1},W_{2},...,W_{n})$\textit{\ with }$W_{i}\in $\textit{int}$%
R_{+}^{k_{i}}\backslash \{0\}$\textit{\ for }$i\in I,$\textit{\ then }$G$%
\textit{\ has at least a Pareto equilibrium point }$\widehat{x}\in X.$

\textit{Proof.} By Theorem 2, $G$ has at least weighted Nash equilibrium
point $\widehat{x}\in K$ with respect of the weighted vector $W.$ Lemma 2
and Remark 3 shows that $\widehat{x}$ is also a weak Pareto equilibrium
point of $G,$ and a Pareto equilibrium point of $G$ if $%
W=(W_{1},W_{2},...,W_{n})$ with $W_{i}\in $int$\mathbb{R}_{+}^{k_{i}}%
\backslash \{0\}$ for each $i\in I.$

\textbf{Acknowledgment:} This work was supported by the strategic grant
POSDRU/89/1.5/S/58852, Project "Postdoctoral programme for training
scientific researchers" cofinanced by the European Social Found within the
Sectorial Operational Program Human Resources Development 2007-2013.

The author thanks to Professor Jo\~{a}o Paulo Costa from the University of
Coimbra for the fruitfull discussions and for the hospitality he proved
during the visit to his departament.


\begin{thebibliography}{99}
\bibitem{Borm} P. Borm, F. Megen, S. Tijs, \textit{A perfectness concept for
multicriteria games.} Math. Meth. Oper. Res. \textbf{49} (1999), 401-412.

\bibitem{Chebbi} S. Chebbi, \textit{Existence of Pareto equilibria for
non-compact constrained multi-criteria games.} J. Appl. An. \textbf{14}
(2008), \textit{2}, 219-226.

\bibitem{Fan2} K. Fan, \textit{A generalization of Tyhonoff's fixed point
theorem.} Math. Ann. \textbf{142} (1961), 305-310.

\bibitem{Fan} K. Fan, \textit{A minimax inequality and applications} in: O.
Shisha (Ed.), Inequalities III, Academic Press, New York, 1972, pp. 103-113.

\bibitem{Kim} W. K. Kim, \textit{Weight Nash equilibria for generalized
multiobjective games.} J. Chungcheong Math. Soc. \textbf{13} (2000), \textit{%
1,} 13-20.

\bibitem{Kim, Ding} W. K. Kim, X. P. Ding, \textit{On generalized weight
Nash equilibria for generalized multiobjective games.} J. Korean Math. Soc. 
\textbf{40} (2003), \textit{5}, 883-899.

\bibitem{Luc} D. T. Luc, E. Sarabi and A. Soubeyran, \textit{Existence of
solutions in variational relation problems without convexity.} J. Math.
Anal. Appl. \textbf{364} (2010), 544-555.

\bibitem{Park1} S. Park, \textit{On generalizations of the KKM principle on
abstract convex spaces.} Nonlinear Anal. Forum \textbf{11} (2006), \textit{1,%
} 67--77.

\bibitem{Park 2} S. Park, \textit{Elements of the KKM theory on abstract
convex spaces.} J. Korean Math. Soc. \textbf{45} (2008), \textit{1,} 1--27.

\bibitem{Park3} S. Park, \textit{Generalizations of the Nash
EquilibriumTheorem in the KKM Theory. }Fixed Point Theory Appl.
doi:10.1155/2010/234706

\bibitem{Park 4} S. Pak, \textit{The KKM principle in abstract convex
spaces: equivalent formulations and applications.} Nonlinear Anal. \textbf{73%
} (2010), 1028-1042.

\bibitem{Park5} S. Park, \textit{Generalizations of the Nash Equilibrium
Theorem in the KKM Theory}. Fixed Point Theory Appl.,
doi:10.1155/2010/234706.

\bibitem{Park} S. Park, \textit{New generalizations of basic theorems in the
KKM theory.} Nonlinear Anal. 74 (2011), 3000-3010.

\bibitem{Voorneveld} M. Voorneveld, S. Grahn, M. Dufwenberg, \textit{Ideal
equilibria in noncooperative multicriteria games.} Math. Meth. Oper. Res. 
\textbf{52} (2000), 65-77.

\bibitem{S. Y. Wang} S. Y. Wang, \textit{Existence of a Pareto equilibrium,}
J. Optim. Theory Appl. \textbf{95} (1997), 373-384.

\bibitem{Yu} H. Yu, \textit{Weak Pareto Equilibria for Multiobjective
Constrained games. }Appl. Math. Let. \textbf{16} (2003), 773-776.

\bibitem{Yu, Yuan} J. Yu, G. X.-Z Yuan, \textit{The study of Pareto
Equilibria for Multiobjective games by fixed point and Ky Fan Minimax
Inequality methods.} Computers Math. Applic. \textbf{\ 35}, (1998), \textit{%
9,} 17-24.

\bibitem{Yuan, Tarafdar} X. Z. Yuan, E. Tarafdar, \textit{Non-compact Pareto
equilibria for multiobjective games.} J. Math. An. Appl. \textbf{204}
(1996),156-163.

\bibitem{Zeleny} M. Zeleny, \textit{Game with multiple payoffs. }Internat.
J. Game Theory \textbf{4} (1976), 179-191.
\end{thebibliography}
\end{document}